\documentclass[11pt]{amsart}
\usepackage{amssymb,amscd,setspace,verbatim}

\newtheorem{theorem}{Theorem}[section]

\theoremstyle{definition}

\newtheorem{example}[theorem]{Example}
\newtheorem{remark}[theorem]{Remark}
\numberwithin{equation}{theorem}

\newenvironment{Macaulay2}
{\medskip\spacing{0.8}\small\verbatim}
{\endverbatim\normalsize\endspacing\medskip}

\def\FF{\mathbb{F}}
\def\NN{\mathbb{N}}

\def\ZZ{\mathbb{Z}}

\def\ge{\geqslant}
\def\le{\leqslant}
\def\bar{\overline}
\def\to{\longrightarrow}
\def\height{\operatorname{height}}

\begin{document}
\title{An algorithm for computing the~integral~closure}

\author{Anurag K. Singh}
\address{Department of Mathematics, University of Utah, 155 South 1400 East, Salt Lake City, UT~84112, USA}\email{singh@math.utah.edu}

\author{Irena Swanson}
\address{Department of Mathematics, Reed College, 3203 SE Woodstock Blvd, Portland, OR 97202}\email{iswanson@reed.edu}

\thanks{A.K.S.~was supported by NSF grant DMS~0600819.}


\begin{abstract}
We present an algorithm for computing the integral closure of a reduced ring that is finitely generated over a finite field.
\end{abstract}
\maketitle

Leonard and Pellikaan~\cite{LP} devised an algorithm for computing the integral closure of weighted rings that are finitely generated over finite fields. Previous algorithms proceed by building successively larger rings between the original ring and its integral closure,~\cite{deJong, Seidenberg70, Seidenberg75, Stolz, Vasconcelos91, Vasconcelos00}. The Leonard-Pellikaan algorithm instead starts with the first approximation being a finitely generated module that contains the integral closure, and successive steps produce submodules containing the integral closure. The weights in~\cite{LP} impose strong restrictions; these weights play a crucial role in all steps of their algorithm. We present a modification of the Leonard-Pellikaan algorithm which works in much greater generality: it computes the integral closure of a reduced ring that is finitely generated over a finite field.

We discuss an implementation of the algorithm in Macaulay~2, and provide comparisons with de~Jong's algorithm~\cite{deJong}.

\section{The algorithm}

Our main result is the following theorem; see Remark~\ref{rem:conductor} for an algorithmic construction of an element $D$ as below when $R$ is a domain, and for techniques for dealing with the more general case of reduced rings.

\begin{theorem}
\label{thm:main}
Let $R$ be a reduced ring that is finitely generated over a computable field of characteristic~$p>0$. Set $\bar{R}$ to be the integral closure of $R$ in its total ring of fractions. Suppose $D$ is a nonzerodivisor in the conductor ideal of $R$, i.e., $D$ is a nonzerodivisor with $D\bar{R}\subseteq R$.

\begin{enumerate}
\item Set $V_0=\frac{1}{D}R$, and inductively define
\[
V_{e+1}=\{f\in V_e\mid f^p\in V_e\}\qquad\text{ for }e\ge0\,.
\]
Then the modules $V_e$ are algorithmically constructible.

\item The descending chain
\[
V_0\supseteq V_1\supseteq V_2\supseteq V_3\supseteq\cdots
\]
stabilizes. If $V_e=V_{e+1}$, then $V_e$ equals $\bar{R}$.
\end{enumerate}
\end{theorem}

The prime characteristic enables us to use the Frobenius or $p$-th power map; this is what makes the modules $V_e$ algorithmically constructible.

\begin{remark}
\label{rem:idealize}
For each integer $e\ge0$, the module $DV_e$ is an ideal of $R$; we set $U_e=DV_e$ and use this notation in the proof of Theorem~\ref{thm:main} as well as in the Macaulay~2 code in the following section. The inductive definition of $V_e$ translates to $U_0=R$ and
\[
U_{e+1}=\{r\in U_e\mid r^p\in D^{p-1}U_e\}\qquad\text{ for }e\ge0\,.
\]
\end{remark}

\begin{proof}[Proof of Theorem~\ref{thm:main}] (1) By Remark~\ref{rem:idealize}, it suffices to establish that the ideals $U_e$ are algorithmically constructible. This follows inductively since
\[\CD
U_{e+1}=U_e\ \cap\ \ker\big(R@>F>>R@>\pi>>R/D^{p-1}U_e\big)\qquad\text{ for }e\ge0\,,
\endCD\]
where $F$ is the Frobenius endomorphism of $R$, and $\pi$ the canonical surjection. 

\medskip

(2) By construction, one has $V_{e+1}\subseteq V_e$ for each $e$. Moreover, it is a straightforward verification that
\[
V_e=\{f\in V_0\mid f^{p^i}\in V_0\text{ for each }i\le e\}\,.
\]
Suppose $f\in\bar{R}$. Then $f^{p^i}\in\bar{R}$ for each $i\ge0$, so $Df^{p^i}\in R$. It follows that $f\in V_e$ for each $e$.

If $V_{e+1}=V_e$ for some positive integer $e$, then it follows from the inductive definition that $V_{e+i}=V_e$ for each $i\ge 1$.

Let $v_1,\dots,v_s\colon R\to\ZZ\cup\{\infty\}$ be the Rees valuations of the ideal $DR$, i.e., $v_i$ are valuations such that for each $n\in\NN$, the integral closure of the ideal $D^nR$ equals
\[
\{r\in R\mid v_i(r)\ge nv_i(D) \text{ for each }i\}\,.
\]
Let $e$ be an integer such that $p^e>v_i(D)$ for each $i$. Suppose $r/D\in V_e$. Then $(r/D)^{p^e}\in V_0$, so $r^{p^e}\in D^{p^e-1}R$. It follows that
\[
p^ev_i(r)\ge\big(p^e-1\big)v_i(D)
\]
for each $i$, and hence that
\[
v_i(r)\ \ge\ v_i(D)-v_i(D)/p^e\ >\ v_i(D)-1
\]
for each $i$. Since $v_i(r)$ is an integer, it follows that $v_i(r)\ge v_i(D)$ for each~$i$, and therefore $r\in\bar{DR}$. But then $r$ belongs to the integral closure of the ideal $D\bar{R}$ in $\bar{R}$. Since principal ideals are integrally closed in $\bar{R}$, it follows that $r\in D\bar{R}$, whence $r/D\in\bar{R}$.
\end{proof}

\begin{remark}
We claim that if $R$ is an integral domain satisfying the Serre condition $S_2$, then each module $V_e$ is $S_2$ as well.

Proceed by induction on $e$. Without loss of generality, assume $R$ is local. Let $x,y$ be part of a system of parameters for $R$. Suppose $yv\in xV_{e+1}$ for an element $v\in V_{e+1}$. Then $yv/x\in V_{e+1}$, i.e., $yv/x\in V_e$ and $y^pv^p/x^p\in V_e$, or equivalently, $yv\in xV_e$ and $y^pv^p\in x^pV_e$. Since $V_e$ is $S_2$ by the inductive hypothesis, it follows that $v\in xV_e$ and $v^p\in x^pV_e$, hence $v\in xV_{e+1}$.
\end{remark}

\begin{remark}
In the notation of Theorem~\ref{thm:main}, suppose $e$ is an integer such that $V_e=V_{e+1}$. We claim that the integral closure of a principal ideal $aR$ is
\[
\{r\in R\mid Dr^{p^i}\,\in\,a^{p^i}R\text{ for each }i\le e+1\}\,.
\]
To see this, suppose $r$ is an element of the ideal displayed above. Then $Dr^p=ga^p$ for some $g\in R$. Since
\[
D(r/a)^{p^i}\in R\qquad\text{ for each }i\le e+1\,,
\]
it follows that
\[
D(g/D)^{p^i}\in R\qquad\text{ for each }i\le e\,.
\]
But then $g/D\in V_e$, which implies that $g/D\in V_i$ for each $i$. Hence $D(r/a)^{p^i}\in R$ for each $i$, equivalently $r\in\bar{aR}$.
\end{remark}

\begin{remark}
\label{rem:conductor}
Let $R$ be a reduced ring that is finitely generated over a perfect field $K$ of prime characteristic~$p$. We describe how to algorithmically obtain a nonzerodivisor $D$ in the conductor ideal of $R$.

\smallskip

\textbf{Case 1.} Suppose $R$ is an integral domain. Consider a presentation of $R$ over $K$, say $R=K[x_1,\dots,x_n]/(f_1,\dots,f_m)$. Set $h=\height(f_1,\dots,f_m)$. Then the determinant of each $h\times h$ submatrix of the Jacobian matrix $(\partial f_i/\partial x_j)$ multiplies $\bar{R}$ into $R$; this may be concluded from the Lipman-Sathaye Theorem (\cite{LS} or \cite[Theorem~12.3.10]{SH}) as discussed in the following paragraph. At least one such determinant has nonzero image in $R$, and can be chosen as the element $D$ in Theorem~\ref{thm:main}. Other approaches to obtaining an element $D$ are via the proof of \cite[Theorem~3.1.3]{SH}, or equivalently, via the results from Stichtenoth's book~\cite{Stichtenoth}.

Let $J$ be the ideal of $R$ generated by the images of the $h\times h$ submatrices of $(\partial f_i/\partial x_j)$. We claim that $J$ is contained in the conductor of $R$. By passing to the algebraic closure, assume $K$ is algebraically closed. After a linear change of coordinates, assume that the $x_i$ are in general position, specifically, that for any $n-h$ element subset $\Lambda$ of $\{x_1,\dots,x_n\}$, the extension $K[\Lambda]\subseteq R$ is a finite integral extension, equivalently that $K[\Lambda]$ is a Noether normalization of $R$. By the Lipman-Sathaye Theorem, the relative Jacobian $J_{R/K[\Lambda]}$ is contained in the conductor ideal. The claim now follows since, as $\Lambda$ varies, the relative Jacobian ideals $J_{R/K[\Lambda]}$ generate the ideal $J$.

\smallskip

\textbf{Case 2.} In the case where $R$ is a reduced equidimensional ring, one may proceed as above and choose $D$ to be the determinant of an $h\times h$ submatrix of $(\partial f_i/\partial x_j)$, and then test to see whether $D$ is a nonzerodivisor. If it turns out that $D$ is a zerodivisor, set
\[
I_1=(0:_R D)\qquad\text{and}\qquad I_2=(0:_R I_1)\,.
\]
Then each of $R/I_1$ and $R/I_2$ is a reduced equidimensional ring, with fewer minimal primes than $R$, and
\[
\bar{R}=\bar{R/I_1}\times\bar{R/I_2}\,.
\]
Hence $\bar{R}$ may be computed by computing the integral closure of each $R/I_i$.

\smallskip

\textbf{Case 3.} If $R$ is a reduced ring that is not necessarily equidimensional, one may compute the minimal primes $P_1,\dots,P_n$ of $R$ using an algorithm for primary decomposition---admittedly an expensive step---and then compute $\bar{R}$ using Case 1 and the fact that
\[
\bar{R}=\bar{R/P_1}\times\dots\times\bar{R/P_n}\,.
\]
\end{remark}

\section{Implementation and examples}

Here is our code in Macaulay~2~\cite{GS}, which uses this algorithm to compute the integral closure.

\medskip\noindent\textbf{Input:} An integral domain $R$ that is finitely generated over a finite field, and, optionally, a nonzero element $D$ of the conductor ideal of $R$.

\medskip\noindent\textbf{Output:} A set of generators for $\bar{R}$ as a module over $R$.

\medskip\noindent\textbf{Macaulay~2 function:}

\begin{Macaulay2}
icFracP = method(Options=>{conductorElement => null})
icFracP Ring := List => o -> (R) -> (
     P := ideal presentation R;
     c := codim P;
     S := ring P;
     if o.conductorElement === null then (
        J := promote(jacobian P,R);
        n := 1;
        det1 := ideal(0_R);
        while det1 == ideal(0_R) do (
           det1 = minors(c,J);
           n = n+1
        );
        D := det1_0;
     ) else D = o.conductorElement;
     p := char(R);
     K := ideal(1_R);
     U := ideal(0_R);
     F := apply(generators R, i-> i^p);
     while (U != K) do (
        U = K;
        L := U*ideal(D^(p-1));
        f := map(R/L,R,F);
        K = intersect(kernel f, U);
     );
     U = mingens U;
     if numColumns U == 0 then {1_R}
     else apply(numColumns U, i-> U_(0,i)/D)
     )
\end{Macaulay2}

Since the Leonard-Pellikaan algorithm uses the Frobenius endomorphism, it is less efficient when the characteristic of the ring is a large prime. In the examples that follow, the computations are performed on a MacBook~Pro computer with a 2~GHz Intel Core Duo processor; the time units are seconds. The comparisons are with de~Jong's algorithm~\cite{deJong} as implemented in the program \texttt{ICfractions} in Macaulay~2, version 1.1.

\begin{example}
Let $\FF_2[x,y,t]$ be a polynomial ring over the field $\FF_2$, and set $R=\FF_2[x,y,x^2t,y^2t]$. Then $R$ has a presentation
\[
\FF_2[x,y,u,v]/(x^2v-y^2u)\,,
\]
which shows, in particular, that $x^2$ is an element of the conductor ideal. Setting $D=x^2$, the algorithm above computes that the integral closure of $R$ is generated, as an $R$-module, by the elements $1$ and $xyt$. Tracing the algorithm, one sees that $V_0$ is not equal to $V_1$, that $V_1$ is not equal to $V_2$, and that $V_2 = V_3$. Indeed, these $R$-modules are
\[
V_0=\frac{1}{x^2}R\,, \qquad V_1=\frac{1}{x}R+ytR\,,\qquad V_e=R+ytR\ \text{ for } e\ge 2\,.
\]

As is to be expected, the algorithm is less efficient as the characteristic of the ground field increases:

\begin{table}[htdp]
\caption{Integral closure of $\FF_p[x,y,u,v]/(x^2v-y^2u)$}
\begin{center}
\begin{tabular}{|l|r|r|r|r|r|r|r|r|r|}
\hline
characteristic $p$&2&3&5&7&11&13&17&37&97\\
\hline
\texttt{icFracP} &0.04&0.03&0.04&0.04&0.04&0.05&0.05&0.13&0.59\\
\hline
\texttt{icFractions} &0.08&0.09&0.09&0.09&0.14&0.15&0.15&0.15&0.15\\
\hline
\end{tabular}
\end{center}
\end{table}

We remark that $R$ is an affine semigroup ring, so its integral closure may also be computed using the program \texttt{normaliz} of Bruns and Koch~\cite{BK}.
\end{example}

\begin{example}
Consider the hypersurface
\[
R=\FF_p[u,v,x,y,z]/(u^2x^4+uvy^4+v^2z^4)\,.
\]
It is readily verified that $R$ is a domain, and that $t=ux^4/v$ is integral over~$R$. The ring $R[t]$ has a presentation
\[
\FF_p[u,v,x,y,z,t]/I\,,
\]
where $I$ is the ideal generated by the $2\times 2$ minors of the matrix 
\[
\begin{pmatrix}u&t&-z^4\\v&x^4&t+y^4\end{pmatrix}\,.
\]
Since the entries of the matrix form a regular sequence in $\FF_p[u,v,x,y,z,t]$, the ring $R[t]$ is Cohen-Macaulay. Moreover, if $p\neq2$, then the singular locus of $R[t]$ is $V(t,y,xz,vz,ux)$ which has codimension $2$, so $R[t]$ is normal.

If $p=2$ then the ring $R[t]$ is not normal; indeed, in this case, the integral closure of $R$ is generated, as an $R$-module, by the elements
\[
1\,,\qquad\sqrt{uv}\,,\qquad\frac{ux+z\sqrt{uv}}{y}\,,\qquad\frac{vz+x\sqrt{uv}}{y}\,,\qquad\frac{uxz+z^2\sqrt{uv}}{uy}\,.
\]
For small values of $p$, these computations may be verified on Macaulay~2 using either algorithm; some computations times are recorded next. Here, and in the next example, $*$ denotes that the computation did not terminate within six hours.

\begin{table}[htdp]
\caption{Integral closure of $\FF_p[u,v,x,y,z]/(u^2x^4+uvy^4+v^2z^4)$}
\begin{center}
\begin{tabular}{|l|r|r|r|r|r|}
\hline
characteristic $p$&2&3&5&7&11\\
\hline
\texttt{icFracP}&0.07&0.22&9.67&143&12543\\
\hline
\texttt{icFractions}&1.16&$*$&$*$&$*$&$*$\\
\hline
\end{tabular}
\end{center}
\end{table}
\end{example}

\begin{example}
Consider the hypersurface
\[
R=\FF_p[u,v,x,y,z]/(u^2x^p+2uvy^p+v^2z^p)\,,
\]
where $p$ is an odd prime. We shall see that $\bar{R}$ has $p+1$ generators as an $R$-module, but first some comparisons:

\begin{table}[htdp]
\caption{Integral closure of $\FF_p[u,v,x,y,z]/(u^2x^p+2uvy^p+v^2z^p)$}
\begin{center}
\begin{tabular}{|l|r|r|r|r|r|r|r|r|}
\hline
characteristic $p$&3&5&7&11&13&17&19&23\\
\hline
\texttt{icFracP}&0.07&0.09&0.27&1.81&4.89&26&56&225\\
\hline
\texttt{icFractions}&1.49&75.00&4009&$*$&$*$&$*$&$*$&$*$\\
\hline
\end{tabular}
\end{center}
\end{table}

We claim that $\bar{R}$ is generated, as an $R$-module, by the elements
\begin{equation}
\label{eqn:gens}
1,\quad\sqrt{y^2-xz}\,,\quad\text{ and }\quad u^{i/p}v^{(p-i)/p}\text{ for }1\le i\le{p-1}\,.
\end{equation}
It is immediate that these elements are integral over $R$; to see that they belong to the fraction field of $R$, note that
\[
\sqrt{y^2-xz}=\pm\frac{uy^p+vz^p}{u(y^2-xz)^{(p-1)/2}}
\]
and that, by the quadratic formula, one also has
\begin{equation}
\label{eqn:root}
\left(\frac{u}{v}\right)^{1/p}=\frac{-y\pm\sqrt{y^2-xz}}{x}\,.
\end{equation}
Moreover, using~\eqref{eqn:root}, it follows that
\[
v^{1/p}\sqrt{y^2-xz}=\pm(xu^{1/p}+yv^{1/p})\,,
\]
and hence the $R$-module generated by the elements~\eqref{eqn:gens} is indeed an $R$-algebra. It remains to verify that the ring
\[
A=R\big[\sqrt{y^2-xz\,},\ u^{i/p}v^{(p-i)/p}\mid 1\le i\le{p-1}\big]
\]
is normal. For this, it suffices to verify that
\[
B=R\big[\sqrt{y^2-xz\,},\ u^{1/p}\,,v^{1/p}\big]
\]
is normal, since $A$ is a direct summand of $B$ as an $A$-module: use the grading on $B$ where $\deg x=\deg y=\deg z=0$ and $\deg u^{1/p}=1=\deg v^{1/p}$, in which case $A$ is the $p$-th Veronese subring $\bigoplus_{i\in\NN}B_{ip}$. The ring $B$ has a presentation $\FF_p[x,y,z,d,s,t]/I$, where $I$ is generated by the $2\times 2$ minors of the matrix
\[
\begin{pmatrix}y+d&z&s\\x&y-d&-t\end{pmatrix}\,,
\]
and $s\mapsto u^{1/p}$, $t\mapsto v^{1/p}$, $d\mapsto\sqrt{y^2-xz}$. But then---after a change of variables---$B$ is a determinantal ring, and hence normal.
\end{example}

\noindent\textbf{Acknowledgment.} We are very grateful to Douglas Leonard for drawing our attention to \cite{LP} and answering several questions, to Wolmer Vasconcelos for his feedback, and to Amelia Taylor for valuable discussions and help with Macaulay~2.



\begin{thebibliography}{19}
\bibitem{BK}
W. Bruns and R. Koch, \emph{Computing the integral closure of an affine semigroup}, Univ. Iagel. Acta Math. \textbf{39} (2001), 59--70.

\bibitem{deJong}
T. de Jong, \emph{An algorithm for computing the integral closure}, J. Symbolic Comput. \textbf{26} (1998), 273--277.

\bibitem{GS}
D. R. Grayson and M. E. Stillman, \emph{Macaulay~2, a software system for research in algebraic geometry}, available at {\tt http://www.math.uiuc.edu/Macaulay2/}.

\bibitem{LP}
D. A. Leonard and R. Pellikaan, \emph{Integral closures and weight functions over finite fields}, Finite Fields Appl. \textbf{9} (2003), 479--504.

\bibitem{LS}
J. Lipman and A. Sathaye, \emph{Jacobian ideals and a theorem of Brian\c con-Skoda}, Michigan Math. J. \textbf{28} (1981), 199--222.

\bibitem{Seidenberg70}
A. Seidenberg, \emph{Construction of the integral closure of a finite integral domain}, Rend. Sem. Mat. Fis. Milano \textbf{40} (1970), 100--120.

\bibitem{Seidenberg75}
A. Seidenberg, \emph{Construction of the integral closure of a finite integral domain. II}, Proc. Amer. Math. Soc. \textbf{52} (1975), 368--372.

\bibitem{Stichtenoth}
H. Stichtenoth, \emph{Algebraic function fields and codes}, Universitext, Springer-Verlag, Berlin, 1993.

\bibitem{Stolz}
G. Stolzenberg, \emph{Constructive normalization of an algebraic variety}, Bull. Amer. Math. Soc. \textbf{74} (1968), 595--599.

\bibitem{SH}
I. Swanson and C. Huneke, \emph{Integral closure of ideals, rings, and modules}, London Math. Soc. Lecture Note Ser. \textbf{336}, Cambridge Univ. Press, Cambridge, 2006.

\bibitem{Vasconcelos91}
W. Vasconcelos, \emph{Computing the integral closure of an affine domain}, Proc. Amer. Math. Soc. \textbf{113} (1991), 633--638.

\bibitem{Vasconcelos00}
W. Vasconcelos, \emph{Divisorial extensions and the computation of integral closures}, J. Symbolic Comput. \textbf{30} (2000), 595--604.
\end{thebibliography}
\end{document}